\documentclass[11pt]{article}
\usepackage{amsmath,amssymb,epsfig}

\newtheorem{conj}{Conjecture}

\begin{document}

\title{A Conjecture about Molecular Dynamics}
\author{P. F. Tupper}
\date{\today}
\maketitle

\begin{abstract}
An open problem in numerical analysis is to explain why molecular dynamics 
works.  The difficulty is that numerical trajectories are only 
accurate for very short times, whereas the simulations are performed over 
long time intervals.  It is believed that statistical information from 
these simulations is accurate, but no one has offered a rigourous proof of 
this.  In order to give mathematicians a clear goal in understanding this 
problem, we state a precise mathematical conjecture about molecular 
dynamics simulation of a particular system.  We believe that if the 
conjecture is proved, we will then understand why molecular 
dynamics works.
\end{abstract}

\section{Introduction}

Molecular dynamics is the computer simulation of a material at the
atomic level.   In principle the only inputs to a simulation are the
characteristics of a set of particles and a description of the forces
between them.  An initial condition is chosen and from these first
principles the evolution of the  system  in time is simulated using
Newton's laws and a simple numerical integrator \cite{frenkel,allen}.  

Molecular dynamics is a very prevalent computational practice, as a
glance at an issue of the Journal of Chemical Physics will show.   It
does have its limitations: the motion of only a relatively small
number of particles can be simulated over a short time interval. 
However, most of the mesoscopic models that have been suggested to
overcome these difficulties still rely on molecular dynamics as a form
of calibration.  It is likely that molecular dynamics will continue to
be important in the future. 

Given its scientific importance there is very little rigourous
justification of molecular dynamics simulation.  From the viewpoint of
numerical analysis it is surprising that it works at all.  The problem
is that individual trajectories computed by molecular dynamics
simulations are accurate for only small time intervals.  As we will
see in Section~\ref{sec:problem}, numerical trajectories diverge
rapidly from true trajectories given the step-lengths used in
practice.  No one disputes this fact, and no one is particularly
concerned with it either.  The reason is that practitioners are never
interested in particular trajectories to begin with.  They are
interested in ensembles of trajectories.  As long as the numerical
trajectories are representative of a particular ensemble of true
trajectories, researchers are content.  However, that this statistical
information is computed accurately has yet to be rigourously
demonstrated in representative cases. 

The goal of this article is to present a concise mathematical
conjecture that encapsulates this fundamental difficulty.  We present
a model system that is representative of systems  
commonly simulated in molecular dynamics.  We present the results of
numerical simulations of this system using the St\"ormer-Verlet
method, the work-horse of molecular dynamics.  In each simulation a
random initial condition is generated, an approximate trajectory for
the system is computed and the net displacement of one particle over
the duration of the simulation is recorded.  We show that even for
step-sizes that are far too large to accurately compute the position
of the particle, the distribution of the particle's displacement over
the many initial conditions appears to be accurate. 	
From the numerical data we conjecture a rate of convergence for this
particular statistical property.  We believe that if this conjectured
rate of convergence (or one like it) can be rigourously established,
even for this single system, then we will understand significantly
better why molecular dynamics works.

The problem of explaining the accuracy of molecular dynamics
simulation is well-known both in the physical sciences (for example
\cite[p. 81]{frenkel}) and in the mathematics community
\cite{SigStu}. 
This latter reference is a survey of the relation between computation
and statistics for initial value problems in general.   
There has been plenty of 
excellent mathematical work that has done much to explain
various features of this type of simulation, but has not resolved the
issue we consider here.  See \cite{SkeTup,LeiRei,LeBris} for surveys.   

One body of work that has addressed the statistical accuracy of
under-resolved trajectories in a special case is by A.~Stuart and
co-workers.  In \cite{Stu1,Stu2} they 
have explored some linear test systems with provable statistical
properties in the limit of large numbers of particles.  They are able
to show that if the systems are simulated with appropriate methods the
statistical features of numerical trajectories are accurate in the
same limit even when the step-lengths are too large to resolve
trajectories. 
Though these results are interesting since they are the only ones of
their kind now known, for the highly nonlinear  problems of practical
molecular dynamics very different arguments will be required. 

One subproblem that has been attacked more successfully is that of the
computation of ergodic averages. These are averages of functions along
very long trajectories.  All that numerical trajectories have to do to
get these correct is sample the entire phase space evenly.  This is a
much weaker property than getting all statistical features correct.
The most striking work on this question is by S.~Reich \cite{Reich}
which establishes rapid convergence of ergodic averages for
Hamiltonian systems which are uniformly hyperbolic on sets of constant
energy.  Unfortunately, this property has never been established for
realistic systems, and is unlikely to hold for them
\cite{Mackay,Liverani}.  The work \cite{Tupper} established similar
results for systems with much weaker properties but requires radically
small time steps for convergence to occur.

The contribution of this work is to precisely specify a simple problem
 which encapsulates all the essential difficulties of the more general
 problem.   In Section~\ref{sec:system} we present the system we will
 study.  Section~\ref{sec:problem} shows the results of some
 numerical experiments on this system.  There we state our conjecture
 based on the results.  In Section~\ref{sec:approaches} we will
 discuss two possible approaches to proving the conjecture. 
 Finally, in Section~\ref{sec:discussion} we will discuss prospects
 for the eventual resolution of the conjecture. 

\section{The System} \label{sec:system}

\begin{figure}
\epsfig{file=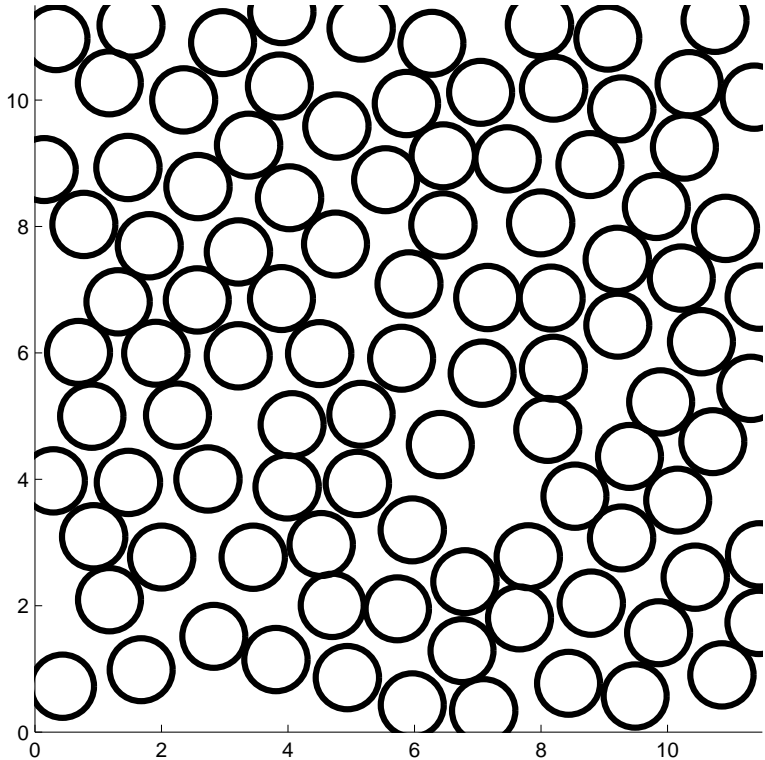,width=4in}
\caption{\label{fig:movie} The positions of the particles for a
  representative state of the system.}   
\end{figure}

The system consists of $n=100$ point particles interacting on an 11.5
by 11.5 square periodic domain.  We let $q \in \mathbb{T}^{2n}$ and $p
\in \mathbb{R}^{2n}$ denote the positions and velocities of the
particles, with $q_i \in \mathbb{T}^2, p_i \in \mathbb{R}^2$ denoting
the position and velocity of particle $i$.   The motion of the system
is described by a system of Hamiltonian differential equations: 
\[
\frac{dq}{dt} = \frac{\partial H}{\partial p}, \ \ \ 
\frac{dp}{dt} = -\frac{\partial H}{\partial q},
\]
with Hamiltonian 
\[
H(q,p) = \frac{1}{2} \|p\|_2^2 + \sum_{i<j} V_{LJ} ( \| q^i-q^j \| ).
\]
Here $V_{LJ}$ denotes the famous Lennard-Jones potential.  In our
simulations we use a truncated version: 
\[
V_{LJ}(r) = \left\{  
\begin{array} {ll}
4 \left( \frac{1}{r^{12}} - \frac{1}{r^6} \right) , & \mbox{if } r
\leq r_{\mbox{\tiny{cutoff}}}, 
\\
0,& \mbox{otherwise.}
\end{array} \right.
\]
Figure~\ref{fig:movie} shows the positions of the particles on the
periodic domain for one state of the system.  Though the particles are
only points, in the figure each is represented by a circle of radius
1/2. 

 We take our initial conditions $q^0,p^0$ to be randomly distributed
 according to the probability density function 
\begin{equation} \label{eq:distrib} 
 Z^{-1} e^{-H(q,p)/k\mathcal{T}},
\end{equation}
where $Z$ is chosen so that the function integrates to one.
This is known as the canonical distribution (or ensemble) for the
system at temperature $\mathcal{T}$.   
There is a simple physical interpretation of this distribution:  if
the system is  weakly connected to another very large system at
temperature $\mathcal{T}$, this is the distribution we will find the
original system in after a long period of time.   
In our units  $k = 1$, 
and we choose $\mathcal{T}=1$.

There are many ways of sampling from the canonical distribution at a
given temperature.  For our experiments we generated initial
conditions using Langevin dynamics.  See  \cite{cances} for an
explanation of this technique and a comparison with other methods.  If
done correctly, the precise method of sampling from the canonical
distribution will have no bearing on the results of the experiments we
will present subsequently. 

The numerical method we use for integrating our system is the
St\"ormer-Verlet scheme.   
Given an initial $q_0, p_0$ and a $\Delta t>0$ it generates a sequence
of states $q_n, p_n, n\geq0$ such that $(q_n,p_n) \approx (q(n\Delta
t),p(n\Delta t))$.  The version of the algorithm we use is  
\begin{eqnarray*}
q_{n+1/2} & = & q_n + p_n \Delta t/2, \\
p_{n+1} & = & p_n - \Delta t \nabla V(q_{n+1/2}), \\
q_{n+1} & = & q_n + p_{n+1} \Delta t/2.
\end{eqnarray*}
This is a second-order explicit method.  It is symplectic, and as a
consequence conserves phase space volume \cite{hairer}.   

Finally we have to decide upon our step-length $\Delta t$.  If $\Delta
t$ is too large the energy of the computed solution will increase
rapidly and explode.  In practice, it is observed that for small
enough step lengths energy remains within a narrow band of the true energy
for very long time intervals. (There is extensive theoretical
justification for this phenomenon, see Section~\ref{subsec:bea}).
Practitioners tend to pick a $\Delta t$ as large as possible while
still maintaining this long-term stability on their time interval of
interest. 
For the system and initial conditions we describe here $\Delta t=
0.01$ yields good approximate energy conservation  on the time
interval $[0,100]$.  For our numerical experiments  we will let
$\Delta t$ take this value and smaller.  (The recommended value in
\cite{frenkel}, a standard reference, for this type of system is
$\Delta t = 0.005$.)

\section{The Problem} \label{sec:problem}

\begin{figure} 
\epsfig{file=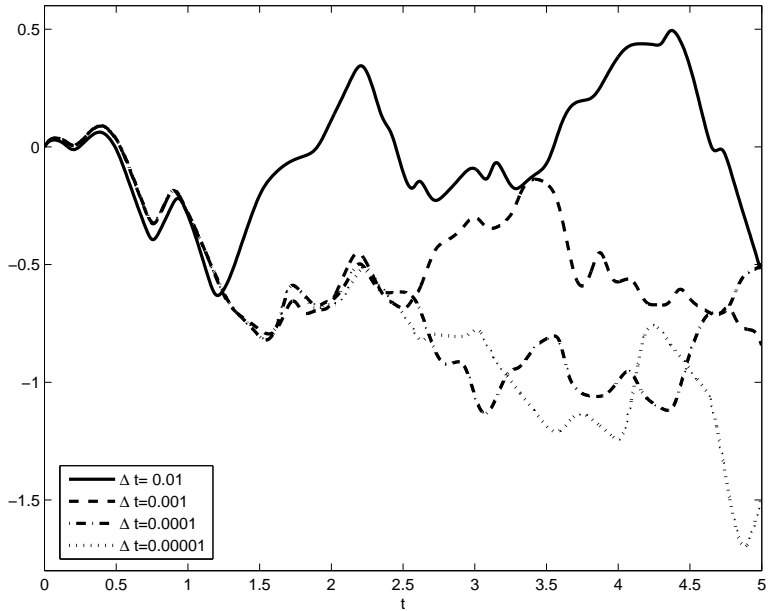,width=4in}
\caption{\label{fig:convergtraj} Computed $x$-position of one particle
  versus time for fixed initial conditions for a  range of $\Delta
  t$.}   
\end{figure}

We will first examine how well trajectories are computed with $\Delta
t = 0.01$. 
Figure \ref{fig:convergtraj} shows  the computed $x$-position of one
particle versus time for the same initial conditions and for a range
of step-lengths.  If the trajectory computed by St\"ormer-Verlet is
accurate over the time interval $[0,5]$, we expect that reducing the
time step by a factor of a thousand would not yield a significantly
different curve.   However, we see that the two curves for $\Delta
t=0.01$ and $\Delta t=0.00001$ very quickly diverge.  They are
distinguishable to the eye almost immediately and completely diverge
around $1.2$ time units. 

Reducing the step length to $\Delta t=0.001$ gives a curve that agrees
with the $\Delta t=0.00001$ line longer, but still diverges around
$2.5$ time units.  Similarly, even with $\Delta t=0.0001$ trajectory
is not accurate over the whole interval depicted.  

From these numerical results, we might conjecture that reducing the
step-length by a constant factor 
only extends the duration for which the simulation is accurate by a
constant amount of time. 
This is consistent with theoretical results about the convergence of
numerical methods for ordinary differential equations. 
What is surprising in this case is that the
time-scale on which the trajectories are valid appears to be miniscule
compared to the time-scale on which computation are actually
performed. 
It seems that the trajectories we compute here with stepsize even as
small as $\Delta t=0.00001$ are not accuarate over the whole interval
$[0,5]$ let alone over considerably longer intervals.

Fortunately we almost never care about what one particular trajectory
is doing in molecular dynamics. 
We only care about statistical features of the trajectories when
initial conditions are selected according to some probability
distribution.  Here we will consider the example of self-diffusion.
Self-diffusion is the diffusion of one particular particle through a
bath of identical particles.  We can imagine somehow marking one
particle at time zero and watching its motion through the system.
This single-particle trajectory will depend on the positions and
velocities of all the particles (including itself) at time zero.
Since these are random, the trajectory of the single particle is
random. 

\begin{figure}
\epsfig{file=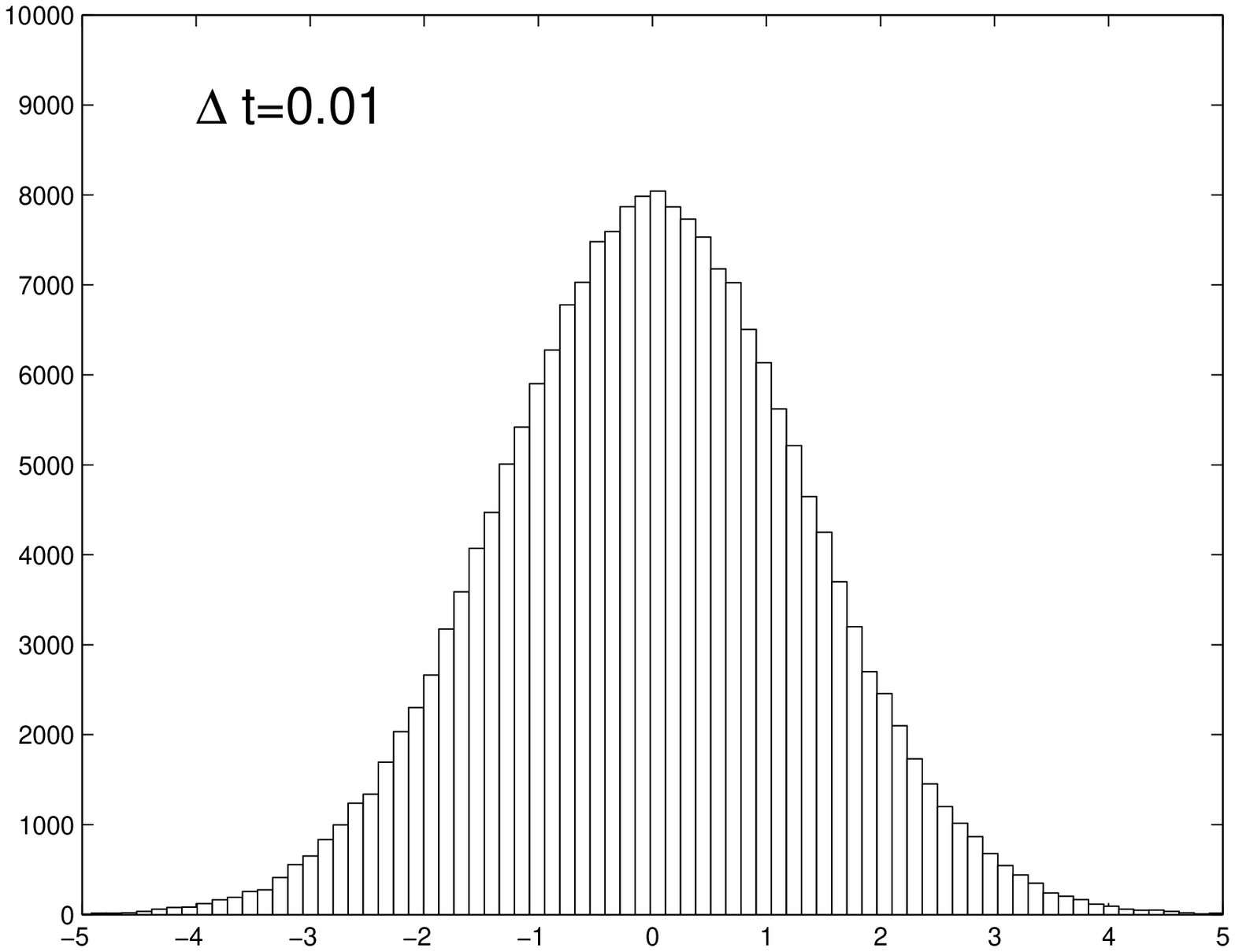,width=1.6in}
\epsfig{file=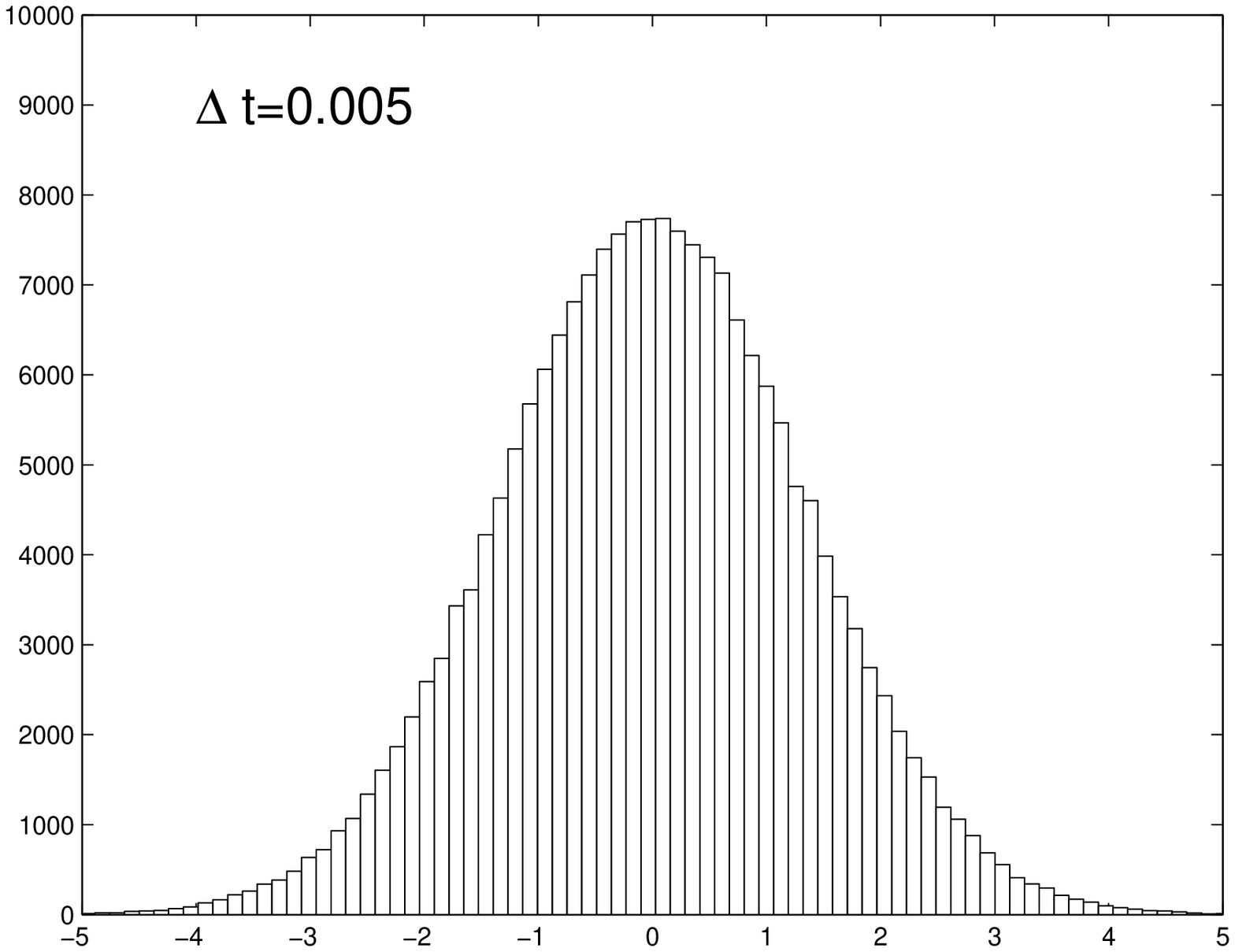,width=1.6in}
\epsfig{file=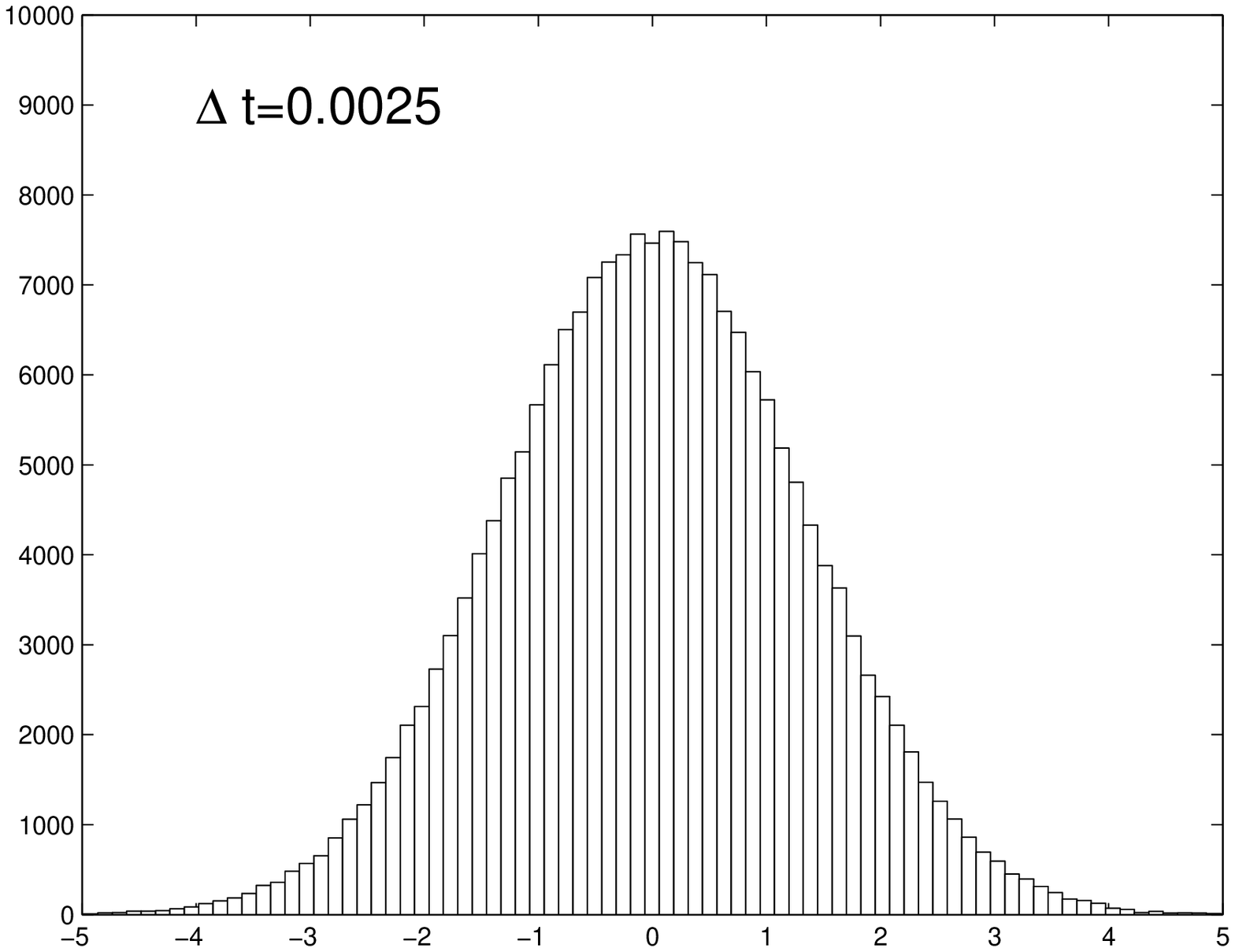,width=1.6in}
\caption{ \label{fig:hists} Displacement in $x$ direction of 1
  particle at $T=10$ for three different step-lengths.} 
\end{figure}

One way to measure self-diffusion is to look at the 
 distribution of the $x$-coordinate  of the tracer particle relative
 to its initial condition.  To estimate this, we generate many random
 initial conditions, perform the simulation using the St\"ormer-Verlet method,
 and record the net displacement of the particle in the given
 direction. 
Figures~\ref{fig:hists} show the histograms of these displacements at
 time $T=10$ for three different step-lengths. 

In contrast to the case where we examined single trajectories, here
the histograms are virtually identical for the different step-lengths.
This suggests that any information we glean from the first histogram
will be accurate. 

To check this more carefully, we compute the variance of the total
displacement at various times $T$ for varying step-lengths.  Let $R(T)
 = \| q_1(T) -  q_1(0) \|$ denote the total displacement of the
particle after time $T$.  This is a random quantity through its
dependence on the state of the system at $t=0$.  Let $R_{\Delta t}(T)$
denote this same displacement as simulated with the St\"ormer-Verlet
method.  This also is a random quantity.  Now define $\langle
R^2_{\Delta t} (T)\rangle$ to be the expected value of $R^2_{\Delta
  t}(T)$ when the initial conditions are chosen according to the
canonical distribution.  Let us see how this last quantity depends on
$\Delta t$.  We do this by generating many initial conditions from the
canonical ensemble and then simulating the system for 100 time units,
keeping track of the total displacement of the tracer particle. 

\begin{figure}
\epsfig{file=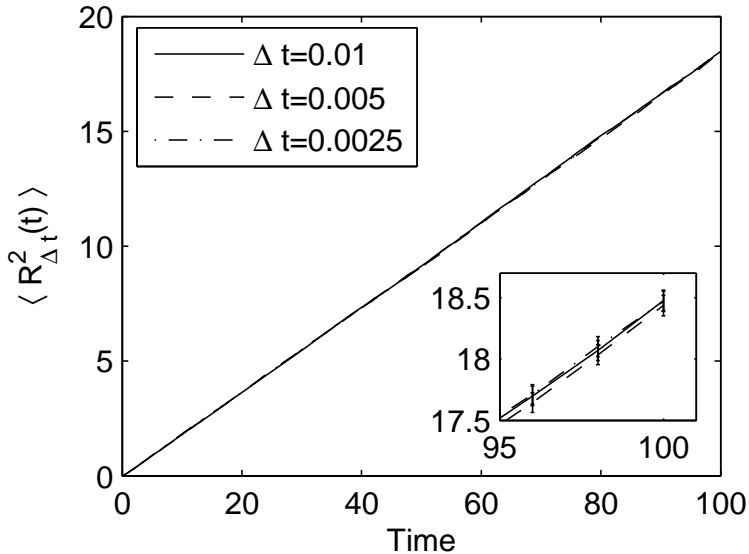,width=4in}
\caption{ \label{fig:selfdiffuse}  Expected squared total displacement in the $x$ direction   of a single particle as a function of time for three different step-lengths. }
\end{figure}

\begin{figure} 
\epsfig{file=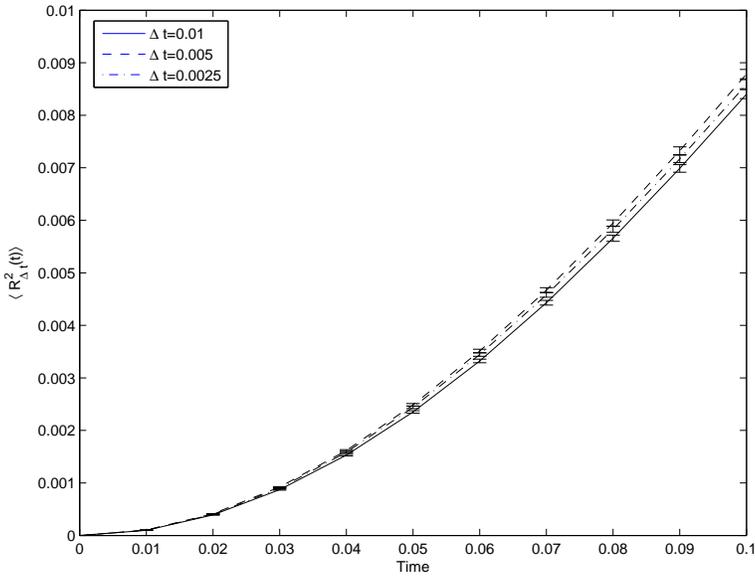,width=4in}
\caption{ \label{fig:zoom} Same as Figure~\ref{fig:selfdiffuse} but on
  a smaller time interval.}
\end{figure}

Figure~\ref{fig:selfdiffuse} shows $\langle R^2_{\Delta t}(T) \rangle$ versus $T$ for three choices of step length.  The inset shows a subset of the data with error bars.  Up to the sampling error there is no difference between the curves.  As far as we can tell from this plot, the answers for $\Delta t=0.01$ are accurate.  The time-scale is much larger than the short interval we found the trajectory to be accurate over.
Lest we give the impression that $\langle R^2_{\Delta t}(T) \rangle$ depends linearly on $T$, Figure~\ref{fig:zoom} shows the same results for a smaller time interval.

We conjecture that the reason $\langle R^2_{\Delta t}(T) \rangle$ does not appear to depend on $\Delta t$ is that even for these large values of $\Delta t$ it closely matches $\langle R^2(T) \rangle$.
It is not clear at all what the rate of convergence of $R_{\Delta t}(T)$ to $R(T)$ is and how it depends on $T$.  However we make the following conjecture:

\begin{conj} \label{conj}
For the system described in Section~\ref{sec:system} with the initial
distribution given by (\ref{eq:distrib}) and the St\"ormer-Verlet
integrator with time step $\Delta t$ 
\[
\left| \langle R^2_{\Delta t}(T) \rangle - \langle R^2(T) \rangle
\right| \leq C \Delta t^2,
\]
for all $T \in [0,A e^{B/\Delta t}]$, for some constants $A, B, C$.
\end{conj}
We will explain the reasons for hypothesizing this particular dependence in the next section.  Here we will briefly note what dependence the classical theory of convergence for numerical ODEs gives:
\[
\left| \langle R^2_{\Delta t}(T) \rangle - \langle R^2(T) \rangle \right| \leq C e^{LT} \Delta t^2
\]
for $T \in [0,E \log(F/\Delta t)]$ for sufficiently small $\Delta t$ for some $C,L,E,F >0$.  (See \cite[p. 239]{stuhum}, for example.)
So we need to explain why the error remains so small even for long
simulations.

\section{Two Approaches} \label{sec:approaches}

We will
discuss two possible approaches to proving Conjecture~\ref{conj}:
backward error analysis and shadowing.

\subsection{Backward Error Analysis} \label{subsec:bea}

Typically a $p$th order numerical method applied to a system of ODEs
computes a trajectory that is $\mathcal{O}(\Delta t^p)$ close to the
exact trajectory on a finite interval. 
Backward error analysis is a way of showing that the numerical
trajectory is an $\mathcal{O}(\exp(-1/\Delta t))$ approximation to the
exact trajectory of a perturbed system.  This result can be used in
turn to prove results about the stability of the numerical trajectory.
See \cite{bennetin} for an early reference and \cite[Ch. IX.]{hairer}
for a recent comprehensive treatment of the subject. 

If we apply a symplectic integrator to a Hamiltonian system it turns
out that the modified system is also Hamiltonian.  The Hamiltonian
function $\widetilde{H}$ for the new system can be written as
$\widetilde{H} = H + \mathcal{O}(\Delta t^2)$.  
There are two consequences for us.  Firstly, the numerical method
agrees very closely with the exact solutions of the modified
Hamiltonian on short time intervals.  If we denote the solution to the
modified system with the same initial conditions by
$(\tilde{q},\tilde{p})$ then 
\begin{equation} \label{eqn:smaller}
| \tilde{q}(n \Delta t) - q^n | \leq C e^{-D/\Delta t}
\end{equation}
for $T \in [0,B/\Delta t]$,
for some appropriate constants \cite{skeel}.  (This alone is not
useful for analysing molecular dynamics since $T$ and $\Delta t$ are
both large.) 
Secondly, the modified Hamiltonian $\widetilde{H}$ is conserved
extremely well by the numerical method for long time intervals: 
\[
\left| \widetilde{H}( q^0,p^0 ) - \widetilde{H}(  q^n,p^n  )  \right|
\leq C e^{-D/\Delta t},
\]
for $n \Delta t \in [0,A e^{B/\Delta t}]$.  Putting this together with
$\widetilde{H} = H + \mathcal{O}(\Delta t^2)$ gives 
\[
\left| H( q^0,p^0 ) - H(  q^n,p^n  )  \right|
\leq E \Delta t^2,
\]
for $n \Delta t \in [0,A e^{B/\Delta t}]$.
We chose the bound in Conjecture 1 in analogy with this last result.

Suppose we wanted to bound the error between $\langle R^2_{\Delta
  t}(t) \rangle$ and $\langle R^2(t) \rangle$ using these estimates.    
The fact that the initial conditions are random adds an extra level of
  complication to the problem.  We have been using $\langle \cdot
  \rangle$ to denote the average with respect to the canonical
  distribution for the Hamiltonian $H$.  The perturbed Hamiltonian
  $\widetilde{H}$ has a different canonical distribution.  We denote
  averages with respect to it by $\langle \cdot \rangle'$.  We let
  $\widetilde{R}$ denote the net displacement of the tracer particle
  under the new flow given by $\widetilde{H}$. 

 We might try bounding the error in the following way:
\begin{eqnarray*}
| \langle R^2_{\Delta t}(T) \rangle - \langle R^2(T) \rangle| & \leq &
| \langle R^2_{\Delta t}(T) \rangle - \langle \widetilde{R}^2(T) \rangle | \\
& & +|\langle \widetilde{R}^2(T) \rangle -  \langle \widetilde{R}^2(T)
 \rangle' | \\ 
& & +|\langle \widetilde{R}^2(T) \rangle' - \langle R^2(T) \rangle|
\end{eqnarray*}
We discuss each of the three terms in turn.

The first term is due to the numerical trajectory not agreeing with
the exact trajectory of the modified system with Hamiltonian
$\widetilde{H}$.  According to (\ref{eqn:smaller}) we can bound this
term by $C \exp(-D/\Delta t)$ for a duration of $B/\Delta t$.  The
studies in \cite{skeel} suggest that this is a tight estimate for
typical molecular dynamics simulations. 

The second term is the difference in the expectation of
$\widetilde{R}^2(t)$ due to a perturbation in the measure.  Since the
two measures are proportional to $\exp(-H/k \mathcal{T})$ and $\exp(
-\widetilde{H}/k \mathcal{T})$ respectively, and $H - \widetilde{H} =
\mathcal{O}(\Delta t^2)$, we expect this term to be on the order of
$\mathcal{O}(\Delta t^2)$ for all $T$.  This probably can be
rigourously controlled without much difficulty.

The third term is just the difference in $\langle R^2(t) \rangle$
between the original system and the perturbed system.  This is likely
to be extremely difficult to bound.  However, showing that it is small
is not a question about computation but about statistical physics.
For now let us assume that 
it is $\mathcal{O}(\Delta t^2)$ for all $T$ for now.

Already we can see that this approach will not get us the result that we
want, even assuming we can bound the third term.
The best estimate we have so far is that
  the error is bounded by $\mathcal{O}(\Delta t^2)$ for $T \in
  [0,B/\Delta t]$.  The bound would hold on an interval much shorter
  than what is needed.
  It appears that backward error analysis alone cannot explain the
  observed convergence.

\subsection{Shadowing} \label{subsec:shadow}

The idea of shadowing is complementary to that of backward error
analysis.  Whereas backward error analysis shows that the numerical
trajectory is close to the exact trajectory of a different Hamiltonian
system with the same initial condition, shadowing attempts to show
that the numerical trajectory is close to an exact trajectory  of the
same Hamiltonian system with a different initial condition.  See
\cite{Hayes} for a nice review of shadowing for Hamiltonian systems. 
  
In our situation, if shadowing were possible, something like the
following would hold. 
Suppose we compute a numerical trajectory starting from $(q^0,p^0)$
with time step $\Delta t$, which we denote by $(q^n,p^n), n\geq 0$.  
 If shadowing is possible then there is an exact trajectory
 $(\tilde{q}(t),\tilde{p}(t))$ of the same Hamiltonian system starting
 at some other initial condition $(\tilde{q}(0),\tilde{p}(0))$ such
 that 
\[
 (q^n,p^n) \approx (\tilde{q}(n \Delta t),\tilde{p}(n \Delta t))
\]
for $n \Delta t$ in some large range of times.
Assuming that it is possible to shadow every numerical trajectory in
this way, let us denote the map on the phase space that takes the
numerical initial condition to the initial condition of the shadow
trajectory by  
\[
S_{\Delta t} (q^0,p^0) = (\tilde{q}(0),\tilde{p}(0)).
\]

The idea of shadowing is used very effectively by S.~Reich in
\cite{Reich}.  For a Hamiltonian system for which shadowing holds he
demonstrates that long-time averages will be computed accurately by
almost all numerical trajectories.  That is,
\begin{equation}  \label{eqn:ergod}
\lim_{T \rightarrow \infty} \frac{1}{T} \int_0^T g(q(t),p(t)) dt
\approx \lim_{N \rightarrow \infty} \frac{1}{N} \sum_{n=0}^N
g(q^n,p^n),
\end{equation}
for almost all initial conditions $(q^0,p^0)=(q(0),p(0))$, for
reasonable functions $g$.  Since the
quantity on the left does not depend on $(q(0),p(0))$ in the systems
considered in \cite{Reich} (except for sets
of measure zero), it is sufficient that
such a map $S_{\Delta t}$ exists to get the result.


In our case we are interested in more general statistical features of
trajectories than long-time averages.
 For example, the variance of
the displacement of a single particle in a finite time interval cannot
be put into the form of a long-time average such as in
$\ref{eqn:ergod}$.  This puts more stringent
requirements on $S_{\Delta t}$. 
To show that statistics are captured correctly 
we cannot consider just single trajectories; we have make sure
the entire ensemble's statistics are reproduced correctly.
If the shadowing map $S_{\Delta t}$ systematically picked initial
conditions for which the tracer particle tended to move to the left,
for example, then the computed statistics could be quite inaccurate.
See \cite{Hayes} for a discussion of this issue in the context of
astrophysics. 
What is necessary for this shadowing to work is for $S_{\Delta t}$ to
leave the canonical ensemble invariant: 
\begin{equation} \label{eqn:preserve}
\langle G(q,p) \rangle = \langle G(S_{\Delta t}(q,p)) \rangle
\end{equation}
for some suitably broad class of functions $G$ on phase space.
This is an even  more stringent requirement than just that shadowing
is possible at all, and it may be quite unlikely to hold for our
system. 

Fortunately we can weaken some other requirements demanded of shadowing 
considerably for our problem.  We do not need the trajectory of the
whole system to be close; we only need the trajectory of a single
particle to be close.  Suppose that our tracer particle's numerical
trajectory is denoted by $(q^n_1,p^n_1)$ for $n\geq 0$.  We say that
\emph{weak shadowing} holds if we can select $\tilde{q}(0)$,
$\tilde{p}(0)$ such that  
\[
 (q^n_1,p^n_1) \approx (\tilde{q}_1(n \Delta t)),\tilde{p}_1(n\Delta t))) 
\]
for $n \Delta t$ in some long range of times.

To see how this fits in with the conjecture suppose that we have both
(\ref{eqn:preserve}) and  
\begin{equation}  \label{eqn:hardthing} 
\| (q^n_1,p^n_1) - (\tilde{q}_1(T)),\tilde{p}_1(T)) \| \leq C \Delta t^2.
\end{equation}
for $T= n \Delta t \in [0, A e^{B/ \Delta t}]$.
This means that  (assuming we can obtain reasonable bounds on
$R_{\Delta t}^2(T)$ and $R^2(T)$) that  
\begin{eqnarray*}
| \langle R^2_{\Delta t} (T) \rangle - \langle R^2(T) \rangle |  
& \leq  & 
K | \langle \| q_1^n \| \rangle - \langle \| q_1(t) \| \rangle |  \\
& \leq &  K | \langle \| q_1^n \| \rangle - \langle \| \tilde{q}_1(T) \| \rangle | +
  K  |  \langle \| \tilde{q}_1(T) \| \rangle      - \langle \| q_1(T) \| \rangle |  \\
& \leq &  K \langle \| (q^n_1,p^n_1) - (\tilde{q}_1(T)),\tilde{p}_1(T)) \| \rangle \\
& & +   K | \langle G(S_{\Delta t}(q^0,p^0)) \rangle - \langle G(q^0,p^0) \rangle |,
\end{eqnarray*}
for $T \in [0,A e^{B/\Delta t}]$.  Here we have let $G$ be the composition of the time $T$ flow map of the Hamiltonian system with the 2-norm.  Now the first term above is bounded by  $C T e^{-D/\Delta t}$ by (\ref{eqn:hardthing}) and the second term is $0$ by (\ref{eqn:preserve}),
thus establishing the conjecture.
Simultaneously proving (\ref{eqn:preserve}) and (\ref{eqn:hardthing}) for some shadowing map $S_{\Delta t}$ may not be easy, but it may be much easier than proving the usual stronger shadowing result.

\section{Discussion} \label{sec:discussion}

Despite the ideas presented in the previous section, the conjecture we
  have presented is probably not open to attack by existing
  techniques.  The problem is that there is no rigourous mathematical
  theory of how statistical regularities emerge from the dynamics of
  generic high-dimensional Hamiltonian systems.  Consequently, there
  is no theory of how perturbations in the Hamiltonian dynamics leads
  to perturbation in the statistics. 
  A numerical analyst has three choices when faced with this situation:
\begin{enumerate}
\item
{\bf Take Up Mathematical Physics.}  
If we are to make progress on the conjecture these entirely
non-numerical problems need to be tackled first.  Mathematical
physicists are interested in proving things like ergodicity and decay
of correlations for Hamiltonian systems such as presented here, and it
is conceivable that eventually there will a robust body of theory that
we can apply to our problem.  So one possibility is to work
on developing such a theory.   This likely will not have much to do
with computation.

\item
{\bf Relax Standards of Rigour.}
Theoretical physicists, as opposed to mathematical physicists, have accepted that much reliable information can be obtained through calculations that cannot be rigourously justified.  Typically theoretical physicists study systems about which nothing interesting can be proved; to do otherwise would be far too restrictive.  There is no reason why this informal yet highly fruitful style of reasoning should be restricted to systems themselves and not numerical discretizations of systems.  A combination of non-rigourous arguments and careful numerical experiments could do a lot to clarify how the St\"ormer-Verlet method is able to compute statistics so accurately for our system.

\item
{\bf Abandon the Whole Pursuit.}  For many, the purpose of numerical analysis is to provide reliable, efficient algorithms.  If one is pursuing a theoretical question, it is hoped that it will lead  to better algorithms eventually.  Sadly, even a complete resolution of the conjecture we have presented in unlikely to have much effect on computational practice.   Many people have tried for years to devise an integrator that is more efficient than the St\"ormer-Verlet method for computing statistically accurate trajectories in molecular dynamics.  They have only been successful for Hamiltonian systems with special structure.  (The prime example of this is the multiple time stepping methods, see \cite[Ch. VIII.4]{hairer}.)
In fact, we state another conjecture which is not formulated rigourously.
\begin{conj}  \label{conj:solver}
No integration scheme can improve the efficiency by more than a factor of two with which St\"ormer-Verlet computes statistically accurate trajectories for systems like that  in Section~\ref{sec:system}.
\end{conj}
Here even a clear mathematical formulation would be a challenge.
Obviously if we already know a lot about a system we can contrive an
algorithm which will give correct statistics for a tracer particle,
but this does not count.  The conjecture is intended to capture the
idea that St\"ormer-Verlet is a very general purpose method; we do not
need to know anything about a system to apply it.\end{enumerate} 

At the Abel Symposium participants seemed to prefer the first of the
three options: try to prove what one can about the system and its
discretization.

{\bf Acknowledgements.} The author was supported by an NSERC Discovery
Grant.  He would like to thank Nilima Nigam, Bob Skeel, and Wayne
Hayes for helpful comments.

\end{document}